# Conjunctive Bayesian networks

NIKO BEERENWINKEL[1], NICHOLAS ERIKSSON[2]
and BERND STURMFELS[3]

[1]*Department of Biosystems Science and Engineering, ETH Zürich, 4058 Basel, Switzerland.*
*E-mail: niko.beerenwinkel@bsse.ethz.ch*

[2]*Department of Statistics, University of Chicago, Chicago, IL 60637, USA.*
*E-mail: eriksson@galton.uchicago.edu*

[3]*Department of Mathematics, University of California, Berkeley, CA 94720-3840, USA.*
*E-mail: bernd@math.berkeley.edu*

Conjunctive Bayesian networks (CBNs) are graphical models that describe the accumulation of events which are constrained in the order of their occurrence. A CBN is given by a partial order on a (finite) set of events. CBNs generalize the oncogenetic tree models of Desper *et al.* by allowing the occurrence of an event to depend on more than one predecessor event. The present paper studies the statistical and algebraic properties of CBNs. We determine the maximum likelihood parameters and present a combinatorial solution to the model selection problem. Our method performs well on two datasets where the events are HIV mutations associated with drug resistance. Concluding with a study of the algebraic properties of CBNs, we show that CBNs are toric varieties after a coordinate transformation and that their ideals possess a quadratic Gröbner basis.

## 1. Introduction

The conjunctive Bayesian network (CBN) model on a finite partially ordered set (poset) was introduced in Beerenwinkel *et al.* ([7], Section 4) as well as in the form of noisy-AND models in the AI literature (e.g., Pearl [20]). Here, we give a self-contained study of the statistical and algebraic properties of this model. CBNs are specializations of Bayesian networks. They include the oncogenetic (also called *mutagenetic*) tree models of Desper *et al.* [10] which have proven very useful in cancer research (Radmacher *et al.* [22]) and in the study of HIV drug resistance (Beerenwinkel *et al.* [5]).

The models are motivated by the following class of problems. Consider a finite set of genetic events, for example, DNA mutations or chromosomal alterations, and assume that the genetic changes are permanent. In this situation, each individual, defined by its genotype, is completely characterized by the subset of the events that have occurred. We wish to learn the constraints on the orders in which these events have accumulated. A CBN is a probabilistic model of this process derived from a partial order on the set of events. This partial order encapsulates the dependencies between events.

For example, consider the development of drug resistance in HIV. This evolutionary process is characterized by the accumulation of resistance mutations in the viral genome. Under fixed drug pressure, these mutations are virtually non-reversible because they confer a strong selective advantage. Thus, the genetic events are fixations of specific amino acid substitutions in the virus population. In each patient, different combinations of resistance mutations will occur. We seek to determine the prevalent mutational pathways along which HIV accumulates drug resistance (cf. Beerenwinkel *et al.* [6]). An order constraint might read that mutations at position 20 and 82 of the target protein must occur before we can see a mutation at position 54. This constraint appears in Figure 4(a). We will analyze such data in Section 4 and see that CBNs are an efficient, accurate tool for this problem.

As another example, the development of cancer is associated with large-scale genetic events such as the gains or losses of parts of chromosomes (Iwasa *et al.* [15]). Knowledge of the constraints on the accumulation of these genetic events helps in assessing the progression of the cancer and assigning treatments (cf. Rahnenführer *et al.* [23]).

A CBN consists of a set of binary random variables, called *events*, and a partial order on these events. While we will use the language of the theory of posets, readers can equivalently think of the partial order as a directed acyclic graph (DAG), with edges encoding the order relations. CBNs are specializations of Bayesian networks, the difference being that in a CBN, an event cannot occur until all of its parents have occurred. Thus, the events that occur with positive probability form a *distributive lattice*. Distributive lattices are important combinatorial objects which have been studied in statistics. For example, the LCI models of (Andersson and Perlman [1] and Andersson *et al.* [2]) use distributive lattices to encode conditional independence statements. Although similar in spirit, readers should beware that LCI models are not the same as CBNs.

Our original motivation for studying CBNs came from work on mutagenetic trees, introduced in Desper *et al.* [10]. Mutagenetic trees assume that each event depends on the occurrence of at most one other previous event. CBNs relax this assumption, allowing for an arbitrary partial order on the events. By relaxing this assumption, CBNs are able to model a larger range of biomedical problems effectively.

Even though they generalize currently used models, CBNs are still very restrictive compared to Bayesian networks in general. However, CBNs have the benefit that the maximum likelihood parameters and structure can be written down in closed form (Proposition 2 and Theorem 5). This is an uncommon phenomenon in the theory of graphical models and should be of independent interest. In addition, the number of parameters in a CBN does not depend on the graph structure, so we do not need to use, for example, the AIC or BIC procedures.

CBNs have also been studied under the name of noisy-AND models in the AI community (Meek and Heckerman [17] and Pearl [20, 21]) as a model for causal inference. The basic idea is that a number of causes influence a common effect through latent intermediate variables; the noisy-AND model requires all causes to have happened before the effect can occur. The study of these models focuses on learning the causal structure given latent variables, in contrast to our situation where we wish to learn the structure of a network of observed variables.

In this paper, we show that CBNs have desirable algebraic, statistical and combinatorial properties. CBN models can be learned efficiently, they can be extended to take into account noise in the data and they perform better than mutagenetic trees in our applications (cf. Figure 3). This paper is organized as follows. After formally introducing CBNs in Section 2, we compute the maximum likelihood (ML) estimator for a CBN in Section 3 and use this to give a combinatorial characterization of the CBN model of maximal likelihood. Next, in Section 4, we compare the performance of CBNs to mutagenetic trees on two data sets of HIV drug resistance mutations. Finally, in Section 5, we study algebraic properties of CBNs. These properties are surprisingly similar to other algebraic results for statistical models. This material may ultimately become relevant for statistical inference, but may also be of independent interest to mathematicians. We determine the prime ideal of algebraic invariants of a CBN and show that this model is a toric variety in a suitable coordinate system. Our main tool is the Möbius transform, a standard tool in working with posets which has found application in the graphical models literature (cf. Drton and Richardson [11], Section 3, and Lauritzen [16], page 239).

## 2. Conjunctive Bayesian networks

A CBN model is specified by a set $\mathcal{E}$ of events, a partial order "$\leq$" on the events and parameters $\theta_e$ for each event $e$. We will assume that there are $n$ events, labeled as $[n] := \{1, \ldots, n\}$. Therefore, we write the parameters as $\theta = (\theta_1, \ldots, \theta_n)$. Frequently, we will abuse notation and refer to both the model $(\mathcal{E}, \leq, \theta)$ and the poset $(\mathcal{E}, \leq)$ as $\mathcal{E}$ when the meaning is clear from the context. A relation $e_1 < e_2$ between two events in $\mathcal{E}$ is interpreted as the requirement that event $e_1$ must happen before event $e_2$ can. The parameter $\theta_e$ is the conditional probability that the event $e \in \mathcal{E}$ has occurred, given that its predecessor events have already occurred.

The state space of the CBN model is the distributive lattice $\mathcal{G} = J(\mathcal{E})$ of order ideals in $\mathcal{E}$. An order ideal is a subset $g \subseteq \mathcal{E}$ such that if $e_2 \in g$ and $e_1 < e_2$, then $e_1 \in g$. Readers unfamiliar with posets and their distributive lattices are referred to Beerenwinkel *et al.* [7], Section 2, for a brief introduction. The elements of $\mathcal{G}$ are called *genotypes*. Thus, a genotype $g \in \mathcal{G}$ is a subset of $\mathcal{E}$ or, equivalently, the binary string in which each bit indicates the occurrence of an event. This terminology presumes a well-defined ground state $0\ldots0$ in which none of the events have yet occurred. In our examples, the unmutated virus strain or the unmutated potential cancer cell is referred to as the "wild type." Hence, for describing mutant types, we need only keep track of which sites differ from the wild type because of the assumption of non-reversibility of events.

We write $\min(g^c)$ for the minimal elements in the complement $g^c = \mathcal{E} \setminus \{g\}$ of a genotype $g$. The elements of $\min(g^c)$ are the events that have not occurred in $g$ but could happen next. For example, in Figure 1, if $g = \{1, 2\}$, then $\min(g^c) = \{3, 4\}$. The probability of observing the genotype $g \in \mathcal{G}$ in the CBN model on the poset $\mathcal{E}$ is defined to

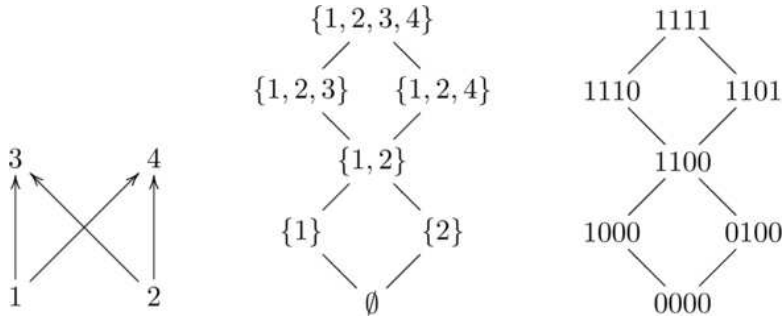

**Figure 1.** Poset on four events, order ideals and genotype lattice.

be

$$P_g(\theta) = \prod_{e \in g} \theta_e \cdot \prod_{e \in \min(g^c)} (1 - \theta_e).$$

That is, the probability of observing $g$ is the probability that all of the events in $g$ have happened times the probability that none of the events that depend only on $g$ have happened.

Equivalently, the CBN model on $\mathcal{E}$ is the directed graphical model for the binary random variables $(X_e)_{e \in \mathcal{E}}$ whose graph has edges $e \to f$ for all cover relations $e \lessdot f$ in $\mathcal{E}$ and whose conditional probability tables are

$$[\Pr(X_e = b \mid X_{\mathrm{pa}(e)} = a)]_{a \in \{0,1\}^{\mathrm{pa}(e)}, b \in \{0,1\}} = \begin{bmatrix} 1 & 0 \\ \vdots & \vdots \\ 1 & 0 \\ 1 - \theta_e & \theta_e \end{bmatrix},$$

where $\mathrm{pa}(e)$ denotes the parents of $e$ in the acyclic directed graph $\mathcal{E}$.

***Example 1.*** Let $n = 4$ and suppose $\mathcal{E}$ is the poset defined on four events by the cover relations $1 < 3$, $1 < 4$, $2 < 3$ and $2 < 4$. The poset $\mathcal{E}$ has precisely seven order ideals, so the distributive lattice $\mathcal{G}$ consists of seven genotypes. They are displayed in Figure 1. The CBN model $\mathcal{E}$ is the family of probability distributions on $\mathcal{G}$ which is given parametrically as follows:

$$P_\varnothing(\theta) = (1 - \theta_1)(1 - \theta_2), \quad P_1(\theta) = \theta_1(1 - \theta_2),$$
$$P_2(\theta) = \theta_2(1 - \theta_1), \quad P_{12}(\theta) = \theta_1 \theta_2 (1 - \theta_3)(1 - \theta_4),$$
$$P_{1234}(\theta) = \theta_1 \theta_2 \theta_3 \theta_4, \quad P_{123}(\theta) = \theta_1 \theta_2 \theta_3 (1 - \theta_4),$$
$$P_{124}(\theta) = \theta_1 \theta_2 \theta_4 (1 - \theta_3).$$

The sum of these seven polynomials equals one. The parameters are the conditional probabilities $\theta_e = \Pr(X_e = 1 \mid X_{\mathrm{pa}(e)} = (1, \ldots, 1))$.

## 3. Maximum likelihood estimation

Consider any CBN model $\mathcal{E}$ on $[n]$. The data for this model take the form of a function $u \colon \mathcal{G} \to \mathbb{N}$, $g \mapsto u_g$, where $u_g$ is the number of observations of the genotype $g$. Given such data $u \in \mathbb{N}^{\mathcal{G}}$, the following proposition gives an easy formula for maximum likelihood estimation of the model parameters.

**Proposition 2.** *For each event $e$ in the CBN model $\mathcal{E}$, the ML estimator $\widehat{\theta}_e$ of $\theta_e$ equals the relative frequency of the genotypes which contain $e$ among all genotypes that contain the events which are strictly below $e$. In symbols,*

$$\widehat{\theta}_e = \frac{\sum_{g: e \in g} u_g}{\sum_{g: \mathrm{below}(e) \subseteq g} u_g} \qquad \textit{for all } e \in \mathcal{E}.$$

**Proof.** The log-likelihood function for the given data $u \in \mathbb{N}^{\mathcal{G}}$ equals

$$\ell_u(\theta) = \sum_{g \in \mathcal{G}} u_g \cdot \left( \sum_{e \in g} \log \theta_e + \sum_{e \in \min(g^c)} \log(1 - \theta_e) \right).$$

The partial derivative of this expression with respect to a parameter $\theta_e$ is

$$\frac{\partial \ell_u}{\partial \theta_e} = \frac{A}{\theta_e} - \frac{B}{1 - \theta_e},$$

where $A$ is the sum over all frequencies $u_g$ of genotypes $g$ containing $e$ and $B$ is the sum over all frequencies $u_g$, where $e \notin g$ but $\mathrm{below}(e) \subseteq g$. Equating this partial derivative with zero, we obtain

$$\widehat{\theta}_e = \frac{A}{A+B},$$

which is precisely the formula asserted in the proposition. □

*Example 3.* We illustrate Proposition 2 for the model in Example 1 and Figure 1. Since $\mathrm{below}(1) = \varnothing$, the ML estimator for $\theta_1$ is

$$\widehat{\theta}_1 = \frac{u_1 + u_{12} + u_{123} + u_{124} + u_{1234}}{u + u_1 + u_2 + u_{12} + u_{123} + u_{124} + u_{1234}}$$

and similarly for $\theta_2$. For $\theta_3$, $\mathrm{below}(3) = \{1, 2\}$ and hence

$$\widehat{\theta}_3 = \frac{u_{123} + u_{1234}}{u_{12} + u_{123} + u_{124} + u_{1234}}.$$

The expression for the ML estimator of $\theta_4$ is similar.

**Remark 4.** Proposition 2 shows that the ML estimator for the CBN model is a rational function of the data. In the language of Catanese *et al.* [8], this says that the *ML degree* of every CBN model is equal to one.

We identify the elements of $\mathcal{G}$ with strings in $\{0,1\}^n$. A probability distribution on $\mathcal{G}$ is thus an element of the $(2^n - 1)$-dimensional simplex $\Delta$ with coordinates indexed by $\{0,1\}^n$. Write $\mathrm{supp}(u)$ for the non-zero coordinates of $u$, that is, for the genotypes that occur in the data set. We say that $u$ *separates the events* if for any two elements $e, f \in [n]$, there exists $g \in \mathrm{supp}(u)$ such that $g \cap \{e, f\}$ is either $\{e\}$ or $\{f\}$. If this is not the case, then we can consider $\{e, f\}$ as a single event and replace $[n]$ by $[n-1]$.

We call any genotype $g \subseteq [n]$ *compatible* with the model $\mathcal{E}$ if $g \in J(\mathcal{E}) = \mathcal{G}$. This is equivalent to $P_g(\theta)$ not being the zero polynomial; see also Beerenwinkel and Drton [3], Definition 14.2. The data $u$ are said to be *compatible* with $\mathcal{E}$ if all $g \in \mathrm{supp}(u)$ are compatible with $\mathcal{E}$. Our next theorem is the main result of this section. It gives a combinatorial solution to the problem of model selection among CBNs. Here, any given data set $u \colon \{0,1\}^n \to \mathbb{N}$ is identified with the corresponding empirical probability distribution in $\Delta$. For such $u \in \Delta$, we can compute the ML estimator $\widehat{\theta}$ for each poset $\mathcal{E}$ on $[n]$. We define the *ML CBN model* for $u$ to be the poset $\mathcal{E}$ for which the log-likelihood $\ell_u(\widehat{\theta})$ has the largest numerical value.

**Theorem 5.** *Let $u \in \Delta$ be a probability distribution which separates the events. There is then a unique largest poset $\mathcal{E}_u$ such that $u$ is compatible with $\mathcal{E}_u$, and the poset $\mathcal{E}_u$ is the unique ML CBN model for $u$.*

Here, "largest poset" refers to the refinement relation among posets on $[n]$, that is, $\mathcal{E} \subset \mathcal{E}'$ means that every relation $e < f$ in $\mathcal{E}$ also holds in $\mathcal{E}'$. Note that this inclusion is reversed for the induced genotype lattices: $\mathcal{E} \subset \mathcal{E}'$ if and only if $\mathcal{G} = J(\mathcal{E}) \supset \mathcal{G}' = J(\mathcal{E}')$.

**Proof of Theorem 5.** The probability $P_g(\theta)$ is identically zero if and only if $g$ is not in $\mathcal{G} = J(\mathcal{E})$. This implies that the likelihood function $\prod_{g \in \mathrm{supp}(u)} P_g(\theta)^{u_g}$ is identically zero if and only if $u$ is not compatible with the poset $\mathcal{E}$. Therefore, we need only consider posets $\mathcal{E}$ such that $u$ is compatible with $\mathcal{E}$.

We claim that there is a unique maximal poset $\mathcal{E}_u$ with which $u$ is compatible. Namely, $\mathcal{E}_u$ is the set of all relations $e_1 < e_2$ such that $g \cap \{e_1, e_2\} \neq \{e_2\}$ for all $g \in \mathrm{supp}(u)$. Note that $\mathcal{E}_u$ is then an antisymmetric relation on $[n]$ because $u$ separates the events. The relation $\mathcal{E}_u$ is transitive because $g \cap \{e_1, e_3\} = \{e_3\}$ implies $g \cap \{e_1, e_2\} = \{e_2\}$ or $g \cap \{e_2, e_3\} = \{e_3\}$. Thus, $\mathcal{E}_u$ is a poset and adding any relation makes $u$ incompatible with it.

It remains to show that if $\mathcal{E}_1 \subset \mathcal{E}_2 \subseteq \mathcal{E}_u$, then $\mathcal{E}_2$ is more likely than $\mathcal{E}_1$. It suffices to show this where $\mathcal{E}_1$ and $\mathcal{E}_2$ differ by only one relation, which we assume without loss of generality to be $1 < 2$. Thus, the events 1 and 2 are incomparable in $\mathcal{E}_1$, but 1 must come before 2 in $\mathcal{E}_2$.

Let $\mathcal{G}_1 = J(\mathcal{E}_1)$ and $\mathcal{G}_2 = J(\mathcal{E}_2)$. We begin by finding the ML parameters for the two models. Write $\widehat{\theta}_e$ for the ML parameters for $\mathcal{G}_1$ and $\widehat{\eta}_e$ for the ML parameters for $\mathcal{G}_2$. According to Proposition 2, we have

$$\widehat{\theta}_e = \frac{\sum_{g \in \mathcal{G}_1 : e \in g} u_g}{\sum_{g \in \mathcal{G}_1 : \text{below}(e) \subseteq g} u_g} \quad \text{and} \quad \widehat{\eta}_e = \frac{\sum_{g \in \mathcal{G}_2 : e \in g} u_g}{\sum_{g \in \mathcal{G}_2 : \text{below}(e) \subseteq g} u_g}.$$

Since $\mathcal{G}_1 \supset \mathcal{G}_2 \supseteq \text{supp}(u)$, the numerators of both expressions are the same, that is, we are summing the counts $u_g$ over all genotypes $g$ that contain $e$.

We claim that $\widehat{\theta}_e = \widehat{\eta}_e$ except when $e = 2$. In both cases, the denominator is the sum over $u_g$ for all genotypes $g$ where $e$ has either already occurred or is allowed to occur next. Since $\mathcal{E}_1$ and $\mathcal{E}_2$ differ in only one relation, $1 < 2$, the denominators are the same (and hence $\widehat{\theta}_e = \widehat{\eta}_e$) unless $e = 2$.

In order to further analyze the ML estimates, we set

$$V_1 = \sum_{g : 1 \in g} u_g, \qquad V_2 = \sum_{g : 2 \in g} u_g,$$

$$N = \sum_{g \in \mathcal{G}_1 : \text{below}(2) \subseteq g} u_g, \qquad M = \sum_{g \in \mathcal{G}_2 : \text{below}(2) \subseteq g} u_g.$$

With this notation, the maximum likelihood parameters are

$$\widehat{\theta}_2 = \frac{V_2}{N} \quad \text{and} \quad \widehat{\eta}_2 = \frac{V_2}{M}.$$

Note that since event 1 always happens in the data before event 2, we have $V_2 \leq V_1$. Since $\mathcal{E}_2$ has more conditions than $\mathcal{E}_1$, we have $M \leq N$ and since event 1 is required to happen before event 2 can in $\mathcal{E}_2$, we have $V_1 \leq M$. Combining these inequalities gives us $V_2 \leq V_1 \leq M \leq N$.

Our analysis will involve the ratios of the ML parameters

$$\frac{\widehat{\theta}_2}{\widehat{\eta}_2} = \frac{M}{N}, \qquad \frac{1 - \widehat{\theta}_2}{1 - \widehat{\eta}_2} = \frac{M}{N} \frac{N - V_2}{M - V_2}. \tag{1}$$

For $i = 1, 2$, the likelihood function for the given distribution $u$ equals

$$L_u(\theta; \mathcal{G}_i) = \prod_g \left( \prod_{e \in g} \theta_e^{u_g} \right) \cdot \left( \prod_{e \in \min(g^c)} (1 - \theta_e)^{u_g} \right).$$

Substitute $\theta_e = \widehat{\theta}_e$ for $i = 1$ and $\theta_e = \widehat{\eta}_e$ for $i = 2$. Our assertion states that

$$L_u(\widehat{\theta}; \mathcal{G}_1) \leq L_u(\widehat{\eta}; \mathcal{G}_2).$$

To prove this, we consider the ratio $L_u(\widehat{\theta}; \mathcal{G}_1) / L_u(\widehat{\eta}; \mathcal{G}_2)$, written as a product over $g \in \text{supp}(u)$, and we examine the four possibilities for $g \cap \{1, 2\}$, given as follows.

Case 1: $g = 00*$. Here, event 2 can happen in $\mathcal{E}_1$, but it cannot yet happen in $\mathcal{E}_2$ since it requires event 1 to happen first. This contributes a factor $(1 - \widehat{\theta}_2)^{u_g}$ to the product over $g$ in $L_u(\widehat{\theta}; \mathcal{G}_1)/L_u(\widehat{\eta}; \mathcal{G}_2)$. Since event 2 has not yet happened, there are no factors $\widehat{\theta}_2/\widehat{\eta}_2$ in this product, so everything else cancels.

Case 2: $g = 01*$. This case cannot happen by compatibility.

Case 3: $g = 10*$. Event 2 has not happened in either case, so all of the terms in the product over $e \in g$ cancel. The same set of events can happen in both $\mathcal{G}_1$ and $\mathcal{G}_2$, so everything in the product over $e \in \min(g^c)$ cancels except the factor with $e = 2$, which occurs both in $L_u(\widehat{\theta}; \mathcal{G}_1)$ and in $L_u(\widehat{\eta}; \mathcal{G}_2)$.

Case 4: $g = 11*$. This case is similar to case 3, except that event 2 has now happened in both cases.

The result of this analysis is the identity

$$\frac{L_u(\widehat{\theta}; \mathcal{G}_1)}{L_u(\widehat{\eta}; \mathcal{G}_2)} = \prod_{g=00*} (1 - \widehat{\theta}_2)^{u_g} \prod_{g=10*} \left(\frac{1 - \widehat{\theta}_2}{1 - \widehat{\eta}_2}\right)^{u_g} \prod_{g=11*} \left(\frac{\widehat{\theta}_2}{\widehat{\eta}_2}\right)^{u_g}. \qquad (2)$$

Note that

$$\sum_{g=10*} u_g + \sum_{g=11*} u_g = V_1 \quad \text{and} \quad \sum_{g=11*} u_g = V_2.$$

Therefore, $\sum_{g=00*} u_g = 1 - V_1$. Substituting (1) into (2), we obtain

$$\frac{L_u(\widehat{\theta}; \mathcal{G}_1)}{L_u(\widehat{\eta}; \mathcal{G}_2)} = \left(\frac{N - V_2}{N}\right)^{1-V_1} \left(\frac{M(N - V_2)}{N(M - V_2)}\right)^{V_1 - V_2} \left(\frac{M}{N}\right)^{V_2} \qquad (3)$$

$$= \frac{M^{V_1}}{N} \cdot \frac{(N - V_2)^{1-V_2}}{(M - V_2)^{V_1 - V_2}}.$$

The following lemma shows that (3) is less than or equal to one for all $0 \leq V_2 \leq V_1 \leq M \leq N \leq 1$. This completes the proof of Theorem 5. $\square$

**Lemma 6.** *If $x, y, a, b$ are real numbers with $0 \leq a \leq b \leq x \leq y \leq 1$, then*

$$\frac{x^b}{y} \cdot \frac{(y - a)^{1-a}}{(x - a)^{b-a}} \leq 1. \qquad (4)$$

**Proof.** We fix $a$ and $b$ and regard the left-hand side of (4) as a function $f_{a,b}(x, y)$ of $x$ and $y$. The two partial derivatives of this function satisfy

$$\frac{\partial f_{a,b}}{\partial x} = \frac{a(x - b)}{x(x - a)} \cdot f_{a,b}(x, y) \quad \text{and} \quad \frac{\partial f_{a,b}}{\partial y} = \frac{a(1 - y)}{y(y - a)} \cdot f_{a,b}(x, y).$$

Both expressions are positive on the triangle $\{(x, y) \in \mathbb{R}^2 : \max(a, b) \leq x \leq y \leq 1\}$, hence $f_{a,b}(x, y)$ is bounded above by $f_{a,b}(1, 1) = (1 - a)^{1-b} \leq 1$. $\square$

We summarize the results of this section in the following algorithm.

***Algorithm 7 (Model selection and parameter estimation for CBN models).***
INPUT: A probability distribution $u \in \Delta$ on the set of genotypes $\{0,1\}^n$.
OUTPUT: The ML CBN model $\mathcal{E}_u$ and the ML parameters $\widehat{\theta}$.
STEP 1: Check whether $u$ separates the $n$ events. If not, group non-distinguished events together, thus decrementing $n$, and replace $u$ by the probability distribution which is induced on the smaller set of genotypes.
STEP 2: Define the poset $\mathcal{E}_u$ on $[n]$ as follows. For any two events $e, f \in [n]$, we set $e < f$ in $\mathcal{E}_u$ if and only if $g \cap \{e, f\} \neq \{f\}$ for all $g \in \mathrm{supp}(u)$.
STEP 3: For each event $e \in [n]$, compute $\widehat{\theta}_e$ by the formula in Proposition 2.
STEP 4: Output the poset $\mathcal{E}_u$ and the vector $\widehat{\theta} \in [0,1]^n$.

## 4. Application to HIV genetic data

The use of Algorithm 7 to obtain the ML CBN model is complicated by the presence of noise in real-world data sets. Any relation $e < f$ between two events $e$ and $f$ will be estimated to be part of the poset $\mathcal{E}$ only if no genotype which contains $f$ but not $e$ has been observed. Thus, the algorithm will miss relations $e < f$ that have strong, but imperfect, support. The problem of noisy data has been analyzed in earlier work on mutagenetic tree models. It can be addressed by explicit error models in an ML framework, as described in Beerenwinkel and Drton ([3], Section 14.2) and Beerenwinkel and Drton [4], Section 3.3. Also, Szabo and Boucher [27] have incorporated an error model directly into the reconstruction algorithm of Desper *et al.* [10].

We propose the following method for constructing a range of CBN models as the error tolerance $\epsilon$ varies. Let $\mathcal{E}_\epsilon$ be the poset on $[n]$ which consists of all relations $e < f$ which are violated by at most a fraction $\epsilon$ of the data. Thus, for $\epsilon = 0$, we recover $\mathcal{E}_u$. Generally, some observations $g \in \mathrm{supp}(u)$ will be incompatible with the model $\mathcal{E}_\epsilon$. These samples are removed prior to ML estimation of the model parameters $\theta$. In order to account for both the compatible and incompatible data, we use a simple mixture model.

Write $\mathcal{G}_\epsilon = J(\mathcal{E}_\epsilon)$ for the genotype space of the model $\mathcal{E}_\epsilon$. We assume that the incompatible genotypes $g \notin \mathcal{G}_\epsilon$ are generated with uniform probability $1/(2^n - |\mathcal{G}_\epsilon|)$. Our mixture model $\mathcal{E}'_\epsilon$ is given parametrically by the event probabilities $\theta_e$ and a mixture parameter $\lambda$ as

$$P'_g(\theta, \lambda) = \begin{cases} \lambda P_g(\theta), & \text{if } g \in \mathcal{G}_\epsilon, \\ (1-\lambda)(2^n - |\mathcal{G}_\epsilon|)^{-1}, & \text{if } g \notin \mathcal{G}_\epsilon \end{cases}$$

for each observation $g \in \{0,1\}^n$. This expression gives an explicit trade-off between a large number of compatible samples and a good model fit.

Since the mixing distributions of the model $\mathcal{E}'_\epsilon$ have disjoint support, the log-likelihood function of the data $u: \{0,1\}^n \to \mathbb{N}$ decomposes as follows:

$$\ell'_u(\theta, \lambda) = \sum_{g \in \mathcal{G}_\epsilon} u_g [\log \lambda + \log P_g(\theta)]$$

$$+ \sum_{g \notin \mathcal{G}_\epsilon} u_g [\log(1-\lambda) - \log(2^n - |\mathcal{G}_\epsilon|)]. \tag{5}$$

**Proposition 8.** *The ML estimators $\widehat{\theta}_e$ of $\theta_e$ under the model $\mathcal{E}'_\epsilon$ are given by Proposition 2. The ML estimator $\widehat{\lambda}$ of $\lambda$ under the model $\mathcal{E}'_\epsilon$ is given by the fraction of the data $u$ which is compatible with the model $\mathcal{E}_\epsilon$. That is,*

$$\widehat{\lambda} = \frac{\sum_{g \in \mathcal{G}_\epsilon} u_g}{\sum_g u_g}.$$

**Proof.** The partial derivatives of (5) with respect to $\theta_e$ are the same as they were in Proposition 2. Next, if we solve

$$0 = \frac{\partial \ell'}{\partial \lambda} = \sum_{g \in \mathcal{G}_\epsilon} \frac{u_g}{\lambda} - \sum_{g \notin \mathcal{G}_\epsilon} \frac{u_g}{1-\lambda},$$

we obtain the above formula for $\widehat{\lambda}$. □

We now apply these methods to mutation data from HIV that was obtained from patients under antiretroviral therapy. The set $\mathcal{E}$ of genetic events consists of seven amino acid alterations in the HIV genome that confer drug resistance. Specifically, as an unordered set,

$$\mathcal{E} = \{\text{K20R, M36I, M46I, I54V, A71V, V82A, I84V}\},$$

where, for example, K20R indicates the amino acid mutation from lysine (K) to arginine (R) at position 20 of the HIV protease. We consider two datasets from the Stanford HIV Drug Resistance Database (Rhee *et al.* [24]), which consist of 112 and 691 observed genotypes under therapy with the protease inhibitors ritonavir (RTV) and indinavir (IDV), respectively.

Previous studies identified correlations and preferred pathways among the resistance mutations (Condra *et al.* [9] and Molla *et al.* [18]). In particular, in Beerenwinkel *et al.* [7], we used mutagenetic trees to infer the underlying dependency structure. The posets are displayed in Figure 2.

For each dataset, we built posets $\mathcal{E}_\epsilon$ for various values of $\epsilon$. For each estimated poset, we report two numbers: the log-likelihood of the data given the mixture model $\mathcal{E}'_\epsilon$ and the mixture parameter $\widehat{\lambda}$ (i.e., the fraction of the data which was explained by the model $\mathcal{E}_\epsilon$). We also calculated these numbers for the mutagenetic trees (Figure 2). These results are shown in Figure 3. Software for building the posets $\mathcal{E}_\epsilon$ and computing the likelihood is available at http://bio.math.berkeley.edu/CBN/.

The two posets that maximize $\ell'_u$ for RTV and IDV, respectively, are displayed in Figure 4. Note that almost all constructed CBNs $\mathcal{E}_\epsilon$ performed better than the mutagenetic trees. In order to estimate the significance of this difference, we repeated the log-likelihood calculation for each poset using 1000 bootstrap samples from the original data. The difference in log-likelihood between these optimal posets and the mutagenetic-tree-induced posets is sufficiently large that their distributions derived from the bootstrap analysis never overlapped. Thus, the difference between the optimal CBN models and the mutagenetic trees is found to be highly significant.

Comparing the optimal CBNs (Figure 4) to the mutagenetic trees (Figure 2) suggests that the mutagenetic trees may induce too many relations and may be handicapped by the requirement that the output is a tree. The posets for RTV share two relations (V82A < M46I and V82A < I54V), while those for IDV share none. The RTV poset [Figure 4(a)] includes the conjunction that both mutations K20R and V82A must occur before I54V, which cannot be represented in a mutagenetic tree. By contrast, the IDV poset [Figure 4(b)] could be represented by a mutagenetic tree, but this tree has not been found by the tree-building procedure of Desper *et al.* [10]. Although the posets and trees do not share many relations, they display a similar structure in that the development of ritonavir resistance is a much more ordered process than for indinavir (see also Beerenwinkel *et al.* [7]).

This comparison suggests several advantages of CBNs. First, they provide better model fits than the posets derived from the mutagenetic tree models. Second, they rely on an ML method both for parameter estimation and for model selection. This stands in contrast to the algorithm of Desper *et al.* [10], which is not an ML procedure. Finally, the perturbed CBNs $\mathcal{E}_\epsilon$ can cover a wide range of fractions of unexplained samples, providing a "parametric" picture of the relations present in the data.

## 5. Algebraic study of the CBN model

In this final, section we study CBNs from the perspective of algebraic statistics. Following Pachter and Sturmfels [19], we regard a CBN as an algebraic variety in a space of dimension $|\mathcal{G}|$. The objective is to compute the prime ideals of all polynomials which vanish on this variety. These polynomials are the *algebraic invariants* of the CBN model.

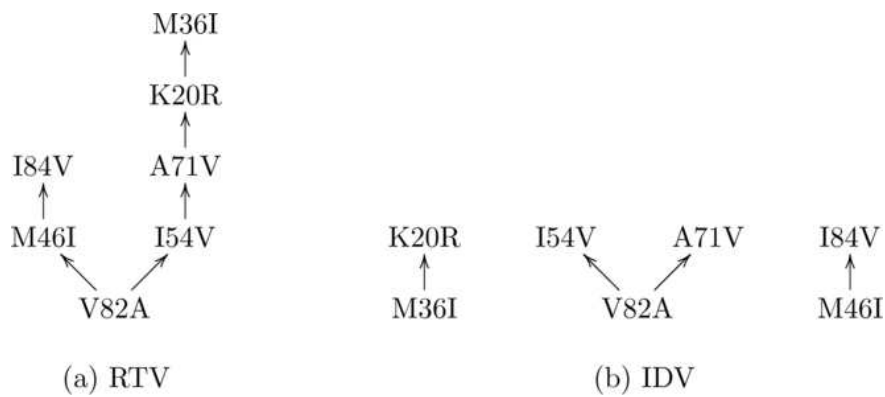

**Figure 2.** Posets corresponding to the mutagenetic trees that were found in Beerenwinkel *et al.* [7], Figure 3, for (a) ritonavir (RTV) and (b) indinavir (IDV).

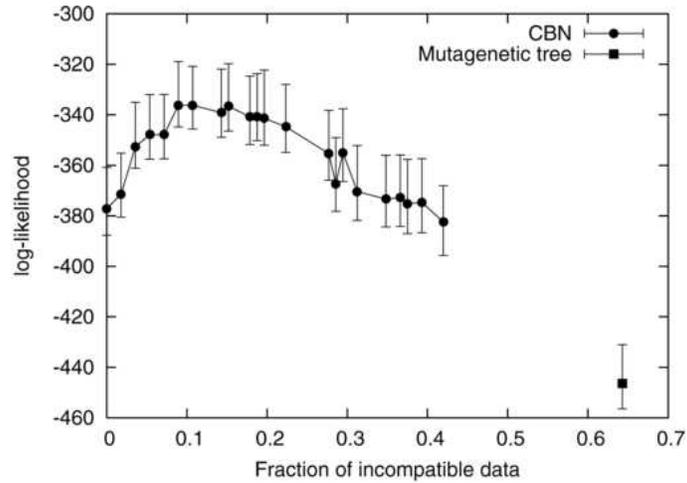

(a) RTV

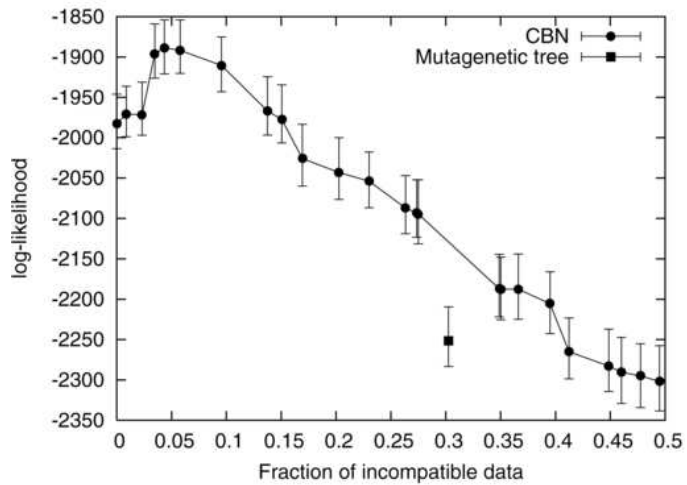

(b) IDV

**Figure 3.** Log-likelihood $\ell'_u$ for the CBN models $\mathcal{E}_\epsilon$ (filled circles) for various choices of the error tolerance $\epsilon$ as a function of the fraction of incompatible genotypes $g \notin \mathcal{G}_\epsilon$. The filled squares correspond to the trees shown in Figure 2. Quartile bars have been derived from 1000 bootstrap samples. Subfigures correspond to (a) ritonavir (RTV) and (b) indinavir (IDV).

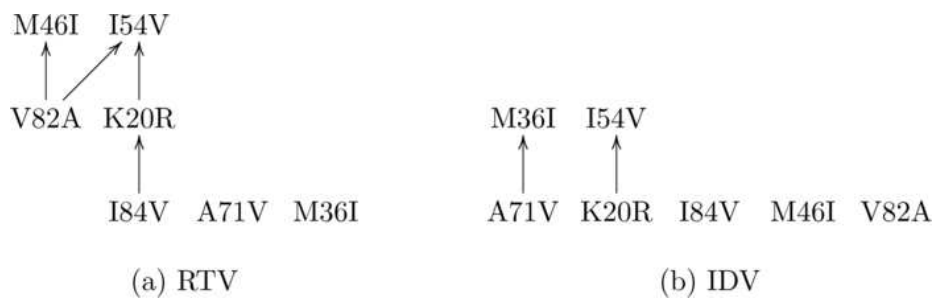

**Figure 4.** Conjunctive Bayesian networks for (a) ritonavir (RTV) and (b) indinavir (IDV) that maximize the likelihood $\ell'_u$. The posets represent the models corresponding to the maxima of the graphs shown in Figure 3.

***Example 9.*** For the model with four events and seven genotypes in Example 1, the algebraic invariants are generated by the three polynomials

$$p_{123} \cdot p_{124} - p_{12} \cdot p_{1234}, \qquad p_1 \cdot p_2 - p_\varnothing \cdot p_{12} - p_\varnothing \cdot p_{123} - p_\varnothing \cdot p_{124} - p_\varnothing \cdot p_{1234}$$

and

$$p_\varnothing + p_1 + p_2 + p_{12} + p_{123} + p_{124} + p_{1234} - 1.$$

Indeed, these three expressions vanish identically if we replace $p_g$ by $P_g(\theta)$ for each genotype $g$. Here "are generated" means that every other polynomial with this property is a linear combination of the three polynomials.

The main theorem in this section exhibits an explicit Gröbner basis for the algebraic invariants of any CBN model. This Gröbner basis consists of a set of quadratic polynomials, together with the trivial invariant $\sum_{g \in \mathcal{G}} p_g - 1$, exactly as in Example 9. For the special case where $\mathcal{E}$ is a forest, this result was proven in Beerenwinkel and Drton [3], Theorem 14.11. Other widely used statistical models have Gröbner bases of the same form, for example, decomposable Markov random fields (Geiger *et al.* [12]) and Jukes–Cantor models in phylogenetics (Sturmfels and Sullivant [26]). Among Markov random fields, having a Gröbner basis of quadrics is equivalent to having ML degree one (Geiger *et al.* [12], Theorem 4.4). This suggests a possible relationship between Theorem 5 and Theorem 10 below. The analogy to Jukes–Cantor models is noteworthy. These models become toric varieties after a linear change of coordinates, known as the *Fourier transform* or *Hadamard conjugation*. The same property will be shown in Corollary 11 for the CBN models, but the role of the Fourier transform is now played by Möbius inversion on the distributive lattice $\mathcal{G}$.

To state our algebraic results, we regard the probabilities $p_g$, for each genotype $g$ in $\mathcal{G} = J(\mathcal{E})$, as unknowns. These generate the polynomial ring

$$\mathbb{R}[\mathcal{G}] = \mathbb{R}[p_g : g \in \mathcal{G}].$$

In this ring, we consider the prime ideal $I_{\mathcal{E}}$ consisting of all polynomials that vanish on the family of probability distributions defined by the CBN model $\mathcal{E}$. Equivalently, $I_{\mathcal{E}}$ is the kernel of the ring map $\mathbb{R}[\mathcal{G}] \to \mathbb{R}[\mathcal{E}], p_g \mapsto P_g(\theta)$, where $\mathbb{R}[\mathcal{E}]$ is the polynomial ring generated by the parameters $\theta_e$, $e \in \mathcal{E}$.

We fix a linear extension of the reverse inclusion order on $\mathcal{G}$, where $g = \varnothing$ is the largest element and $g = \mathcal{E}$ is the smallest element. We define $\prec$ to be the degree reverse lexicographic monomial ordering on $\mathbb{R}[\mathcal{G}]$ induced by the variable ordering given by the fixed linear extension.

**Theorem 10.** *The reduced Gröbner basis of the ideal $I_{\mathcal{E}}$ with respect to the monomial ordering $\prec$ consists of the trivial linear invariant $\sum_{g \in \mathcal{G}} p_g - 1$, with leading term $p_\varnothing$, together with one homogeneous quadratic polynomial*

$$\underline{p_g \cdot p_h} - p_{g \cup h} \cdot p_{g \cap h} + \prec\text{-lower terms} \tag{6}$$

*for each incomparable pair of genotypes $\{g, h\}$ in the distributive lattice $\mathcal{G}$.*

**Proof.** We start our proof of Theorem 10 by recalling that the sum of the polynomials $P_g(\theta)$, equals one. If we take the subsum of all polynomials $P_g(\theta)$, where $g$ runs over all genotypes containing some fixed genotype $h \in \mathcal{G}$, then this is essentially the same sum with $\mathcal{E}$ replaced by $\mathcal{E} \setminus h$ and we conclude that

$$\sum_{g: h \subseteq g} P_g(\theta) = \prod_{e \in h} \theta_e.$$

This expression represents the probability that each event in $h$ has happened. The identity suggests that we perform the following linear change of coordinates in the polynomial ring $\mathbb{R}[\mathcal{G}]$:

$$q_h := \sum_{g: h \subseteq g} p_g \qquad \text{for all } h \in \mathcal{G}. \tag{7}$$

Thus, in the new coordinates $q_h$, the CBN model is precisely the toric variety associated with the distributive lattice $\mathcal{G}$. A well-known theorem of Hibi [14] states that the ideal of this toric variety is generated by the binomials

$$\underline{q_g \cdot q_h} - q_{g \cup h} \cdot q_{g \cap h}, \tag{8}$$

where $\{g, h\}$ runs over all incomparable pairs of elements of $\mathcal{G}$. Moreover, these binomials form a Gröbner basis with the underlined terms as the leading terms. Thus, $I_{\mathcal{E}}$ is generated by the quadrics (8) together with the relation $q_\varnothing - 1$, which is obtained from (7) under the assumption that the probabilities $p_g$ sum to one. Now, if we rewrite (8) in terms of the original coordinates $p_g$, then we obtain quadrics of the form (6).

We claim that the quadrics (8) in the original coordinates $p_g$ form a Gröbner basis for $I_{\mathcal{E}}$. This Gröbner basis will be minimal, but not reduced. We shall verify the Gröbner basis property by using the theory of sagbi bases (or canonical bases), as described by Sturmfels [25], Section 11.

Let $<$ denote a negative degree monomial ordering on the polynomial ring of parameters, $\mathbb{R}[\mathcal{E}] = \mathbb{R}[\theta_e : e \in \mathcal{E}]$. Thus, $<$ is a local monomial ordering in which $1 = \theta^0$ is the largest monomial and monomials of higher total degree are $<$-smaller than monomials of lower total degree. See Greuel and Pfister [13] for an introduction to local monomial orderings.

We shall prove that the coordinate polynomials $P_g(\theta)$ of the CBN model form a sagbi basis for the local ordering $<$, that is, the $<$-leading monomials

$$\text{in}_<(P_g(\theta)) = \prod_{e \in g} \theta_e \qquad (9)$$

generate the algebra of all $<$-leading monomials of polynomials in the image of our ring map $\mathbb{R}[\mathcal{G}] \to \mathbb{R}[\mathcal{E}]$, $p_g \mapsto P_g(\theta)$. Let $J_\mathcal{E}$ be the prime ideal in $\mathbb{R}[\mathcal{G}]$ consisting of all algebraic relations on the initial monomials (9). By Hibi's result, $J_\mathcal{E}$ is generated by the binomials $p_g \cdot p_h - p_{g \cup h} \cdot p_{g \cap h}$ and these binomials form the reduced Gröbner basis of $J_\mathcal{E}$ with respect to $\prec$.

Let $w \in \mathbb{R}^\mathcal{E}$ be a weight vector which represents the local ordering $<$ for the coordinate polynomials $P_g(\theta)$ and let $\mathcal{A}^T w$ be the induced weight vector in $\mathbb{R}^\mathcal{G}$. By Sturmfels [25], Lemma 11.2, we have

$$\text{in}_{\mathcal{A}^T w}(I_\mathcal{E}) \subseteq J_\mathcal{E}. \qquad (10)$$

Importantly, $p_g \cdot p_h - p_{g \cup h} \cdot p_{g \cap h}$ is the initial form of (8) with respect to $\mathcal{A}^T w$, so the reverse inclusion also holds. Thus, equality holds in (10) and the desired sagbi basis property holds by Sturmfels [25], Theorem 11.4.

By Sturmfels ([25], Corollary 11.6(a)), we conclude that the quadratic model invariants (8) form a Gröbner basis of $I_\mathcal{E}$ with respect to $\prec$. This Gröbner basis is minimal and it can be transformed into the reduced Gröbner basis by autoreduction. This completes the proof of Theorem 10. □

A few remarks are in order. The linear transformation between the $p$-coordinates and the $q$-coordinates on the polynomial ring $\mathbb{R}[\mathcal{G}]$, given in (7), is precisely the *Möbius inversion* on the genotype lattice $\mathcal{G}$. Equivalently, the coefficients ($+1$, $-1$ or $0$) of the monomials in the expanded model coordinates $P_g(\theta)$ are precisely the values of the Möbius function on $\mathcal{G}$.

Example 9 illustrates Theorem 10 for the model in Example 1. The three model invariants form a Gröbner basis with leading terms $p_{123} \cdot p_{124}$, $p_1 \cdot p_2$ and $p_\varnothing$. Möbius inversion on the genotype lattice pictured in Figure 1 gives

$$p_\varnothing = q_\varnothing - q_1 - q_2 + q_{12}, \qquad p_1 = q_1 - q_{12}, \qquad p_2 = q_2 - q_{12},$$

$$p_{12} = q_{12} - q_{123} - q_{124} + q_{1234}, \qquad p_{1234} = q_{1234},$$

$$p_{123} = q_{123} - q_{1234}, \qquad p_{124} = q_{124} - q_{1234}.$$

If we perform these substitutions in the reduced Gröbner basis listed in Example 9, then the three given model invariants simplify to

$$q_\varnothing - 1, \qquad q_1 \cdot q_2 - q_\varnothing \cdot q_{12}, \qquad q_{123} \cdot q_{124} - q_{1234} \cdot q_{12}.$$

**Corollary 11.** *The Möbius inversion (7) on the distributive lattice $\mathcal{G} = J(\mathcal{E})$ is a linear change of coordinates which identifies the CBN model $\mathcal{E}$ with the toric variety of the distributive lattice $\mathcal{G}$ defined by Hibi.*

We close with the remark that the sagbi basis property of the coordinate polynomials of the CBN model, which was established in the course of proving Theorem 10, can be used to express any polynomial in the coordinate subalgebra rapidly in terms of the generators $P_g(\theta)$. This process, which is known as *subduction* (Sturmfels [25], Algorithm 11.1), generalizes the classical procedure of expressing any symmetric polynomial in terms of elementary symmetric functions and may be of interest to statisticians.

## Acknowledgements

Part of this work was done while Niko Beerenwinkel was at UC Berkeley and supported by the DFG (BE 3217/1-1). Nicholas Eriksson was supported by the NSF (DMS-06-03448) and Bernd Sturmfels was partially supported by the NSF (DMS-04-56960) and by DARPA (HR0011-05-1-0057).